\documentclass[oneside,12pt]{amsart}

\newtheorem{theorem}{Theorem}
\theoremstyle{definition}
\newtheorem{remark}{Remark}
\newtheorem*{acknowledgments}{Acknowledgments}

\newcommand{\A}{\mathbb{A}}
\newcommand{\C}{\mathbb{C}}
\newcommand{\D}{\mathbb{D}}
\newcommand{\E}{\mathbb{E}}
\newcommand{\R}{\mathbb{R}}
\newcommand{\Z}{\mathbb{Z}}
\renewcommand{\H}{\mathcal{H}}
\renewcommand{\S}{\mathbb{S}}
\newcommand{\Ind}{\operatorname{Ind}}
\newcommand{\Ker}{\operatorname{Ker}}
\newcommand{\p}{\partial}
\newcommand{\II}{\operatorname{II}}

\usepackage{color}

\newcommand{\qmbox}[1]{\quad\mbox{#1}\quad}

\usepackage[margin=3.5cm]{geometry}

\title{Index estimates for planar domains with Robin boundary condition}
\author{Abraão Mendes}
\address{Instituto de Matemática, Universidade Federal de Alagoas,
Maceió-AL, Brazil}

\email{abraao.mendes@im.ufal.br}

\usepackage{pgfplots}
\pgfplotsset{compat=1.18}

\usepackage[T1]{fontenc}
\usepackage{lmodern}

\begin{document}

\begin{abstract}
In this note, we study the index of the Laplace operator on planar domains $\Omega$ with compact smooth boundary, allowing both the compact and noncompact cases, under a natural Robin boundary condition involving the geodesic curvature of $\p\Omega$. We obtain lower bounds for the index in terms of the number of boundary components. Our approach combines conformal and spectral techniques with the topology of the domain, encoded in the space of $L^2$ harmonic vector fields tangent to the boundary. The result is motivated by index estimates for free boundary minimal surfaces and provides an intrinsic counterpart to a borderline geometric situation arising in that setting.
\end{abstract}

\maketitle

\section{Introduction}

In this note, we are interested in estimating the index of the Laplace operator of a not necessarily compact planar domain $\Omega\subset\R^2$ with compact smooth boundary $\p\Omega$ under a natural Robin boundary condition. More precisely, we consider the eigenvalue problem
\begin{align}\label{problema}
\begin{cases}
\Delta u+\lambda u=0 &\mbox{in}\quad\Omega,\\
\p_\nu u+k_gu=0 &\mbox{on}\quad\p\Omega,
\end{cases}
\end{align}
where $\nu$ is the outward unit normal to $\p\Omega$ and $k_g=\langle\nabla_\tau\nu,\tau\rangle$ is the geodesic curvature of $\p\Omega$, with $\tau$ denoting the unit tangent vector to $\p\Omega$.

The quadratic form associated with \eqref{problema} is given by
\begin{align*}
Q(u,u)=\int_\Omega|\nabla u|^2+\int_{\p\Omega}k_gu^2,\quad u\in C_c^\infty(\Omega).
\end{align*}
The \textsl{index} of $\Omega$, denoted by $\Ind(\Omega)$, is the supremum of the dimensions of all linear subspaces $V\subset C_c^\infty(\Omega)$ on which $Q$ is negative definite. If $\Ind(\Omega)=0$, we say that $\Omega$ is \textsl{stable}.

When $\Omega$ is compact, classical spectral theory implies that the eigenvalues of \eqref{problema} form a discrete sequence
\begin{align*}
\lambda_1<\lambda_2\le\cdots\le\lambda_k\nearrow+\infty,
\end{align*}
where $\lambda_1$ is simple and is the only eigenvalue whose eigenfunctions do not change sign.

When $\Omega$ is noncompact, the spectrum of $\Delta$ need not be discrete. Nevertheless, if $\Ind(\Omega)=n$, then problem \eqref{problema} admits $n$ negative eigenvalues, counted with multiplicity, whose associated eigenfunctions belong to $L^2(\Omega)$. We refer the reader to \cite[Proposition~4.9]{HongSaturnino} for a proof of this result in a more general setting.

The motivation for the above problem stems from a borderline geometric\linebreak situation, which we now describe.

Let $M$ be a domain in $\R^3$ with boundary $\p M$, and let $\Sigma$ be an oriented immersed minimal surface in $M$ whose boundary $\p\Sigma$ is compact and meets a smooth portion of $\p M$ orthogonally. We also assume that $\Sigma$ is complete and properly immersed ($\Sigma\cap\p M=\p\Sigma$). In summary, $\Sigma$ is a free boundary minimal surface in $M$.

The \textsl{Morse index} of $\Sigma$ is defined as the index of the quadratic form
\begin{align*}
Q_\Sigma(u,u)
=\int_\Sigma(|\nabla u|^2+|dN|^2u^2)
-\int_{\p\Sigma}\II^{\p M}(N,N)u^2,
\quad u\in C_c^\infty(\Sigma),
\end{align*}
where $\II^{\p M}$ denotes the second fundamental form of $\p M$ in $\R^3$, and $N:\Sigma\to\S^2$ is the Gauss map of $\Sigma$.

Estimating the index of $Q_\Sigma$ from below in terms of the topology of $\Sigma$ is a well-established problem in the theory of minimal surfaces that has been studied by many authors in the settings of both complete noncompact minimal surfaces without boundary and compact free boundary minimal surfaces.

Although we do not aim to provide an exhaustive account of the literature on this problem, we would like to mention the works of Ros \cite{Ros} and Chodosh and Maximo \cite{ChodoshMaximo-I} in the complete setting, Ambrozio, Carlotto, and Sharp \cite{AmbrozioCarlottoSharp-FB} in the compact free boundary setting (in dimension three and higher), and Hong and Saturnino \cite{HongSaturnino} in the capillary setting. In \cite{CavalcanteMendesdosSantos}, the present author, together with Cavalcante and dos Santos, extended the three-dimensional theorem of Ambrozio, Carlotto, and Sharp to the class of noncompact free boundary minimal surfaces with compact boundary. Most of the techniques employed in the present note originate from these works and will therefore be used throughout without further reference.

Assuming that the mean curvature $H^{\p M}$ of $\p M$ is nonnegative, Ambrozio,\linebreak Carlotto, and Sharp were able to show that
\begin{align}\label{eq.AmbrozioCarlottoSharp}
\Ind(Q_\Sigma)\ge\frac{2g+r-1}{3},
\end{align}
where $g$ denotes the genus of $\Sigma$ and $r$ the number of connected components of $\p\Sigma$, in the compact free boundary setting. The proof of \eqref{eq.AmbrozioCarlottoSharp} proceeds by contradiction. Indeed, assuming that
\begin{align*}
3\Ind(Q_\Sigma)<2g+r-1,
\end{align*}
they establish the existence of a nontrivial tangent harmonic one-form $\omega$, whose coordinate functions are then used to derive a highly rigid contradictory situation. For instance, under this assumption they show that $\Sigma$ is totally geodesic and that $H^{\p M}$ vanishes along $\p\Sigma$. Consequently, \eqref{eq.AmbrozioCarlottoSharp} follows whenever either $\Sigma$ is not totally geodesic or $H^{\p M}$ does not vanish identically along $\p\Sigma$. However, the borderline case in which $\Sigma$ is totally geodesic and $H^{\p M}$ vanishes along $\p\Sigma$ appears to require a separate analysis.

If $\Sigma$ is a free boundary totally geodesic surface and $H^{\p M}=0$ along $\Sigma$, then the index estimate reduces to the intrinsic problem \eqref{problema}, since $k_g=-\II^{\p M}(N,N)$ in this setting. In this special situation, we were able both to improve the lower bound in \eqref{eq.AmbrozioCarlottoSharp} and to extend the result to the noncompact setting.

The aim of this note is to prove the following result.

\begin{theorem}\label{thm.main}
Let $\Omega$ be a complete, possibly noncompact domain in $\R^2$ with compact smooth boundary. Let $\Ind(\Omega)$ denote the index of the Laplace operator on $\Omega$ with the Robin boundary condition
\begin{align*}
Bu=\p_\nu u+k_gu=0.
\end{align*}
Then one of the following holds:
\begin{enumerate}
\item $\Omega$ is stable and homeomorphic to a compact disk;
\item $\Omega$ is unstable and
\begin{align*}
\Ind(\Omega)\ge \left\lceil \frac{r-2}{2}\right\rceil,
\end{align*}
where $r$ denotes the number of connected components of $\partial\Omega$.
\end{enumerate}
\end{theorem}

Here $\lceil x\rceil$ denotes the smallest integer greater than or equal to $x$.

In this note, rather than working with harmonic one-forms, we work with the corresponding dual harmonic vector fields.

A smooth vector field $\xi:\Omega\to\R^2$, $\xi=(a,b)$, is said to be
\textsl{harmonic} if
\begin{align*}
\begin{cases}
\operatorname{div}(\xi)=a_x+b_y=0,\\
\langle\nabla_{\p_x}\xi,\p_y\rangle-\langle\nabla_{\p_y}\xi,\p_x\rangle=b_x-a_y=0.
\end{cases}
\end{align*}
In particular, $\Delta\xi=(\Delta a,\Delta b)$ vanishes identically on $\Omega$. We say that $\xi$ is \textsl{tangent to the boundary} $\p\Omega$ if $\langle\xi,\nu\rangle=0$ along $\p\Omega$.

We denote by $\H_T^1(\Omega)$ the vector space of harmonic vector fields on $\Omega$ that are tangent to the boundary, and by $L^2\H_T^1(\Omega)$ the subspace consisting of those $\xi\in\H_T^1(\Omega)$ such that $|\xi|\in L^2(\Omega)$.

It is well known that the dimension of $L^2\H_T^1(\Omega)$ is closely related to the topology of $\Omega$. In fact, one can prove that
\begin{align*}
\dim L^2\H_T^1(\Omega)=r-1.
\end{align*}

\section{Proof of Theorem~\ref{thm.main}}

We divide the proof into four cases.

\medskip

\textbf{Case I: $\Omega$ is compact stable.} In this case, there exists a positive eigenfunction $\phi\in C^\infty(\Omega)$ of $\Delta$ corresponding to the first eigenvalue $\lambda_1\ge0$:
\begin{align*}
\begin{cases}
\Delta\phi+\lambda_1\phi=0&\mbox{in}\quad\Omega,\\
\p_\nu\phi+k_g\phi=0&\mbox{on}\quad\p\Omega.
\end{cases}
\end{align*}

The Gaussian and the geodesic curvatures of $\bar{g}=\phi^2\delta$ satisfy
\begin{align*}
\bar{K}\phi^2=-\Delta\log\phi=\lambda_1+\frac{|\nabla\phi|^2}{\phi^2}\ge0,\quad
\bar{k}_g\phi=k_g+\p_\nu\log\phi=0.
\end{align*}

By the Gauss-Bonnet theorem,
\begin{align*}
0\le\int_\Omega\bar{K}+\int_{\p\Omega}\bar{k}_g=2\pi(2-r),
\end{align*}
which yields $r=1,2$. If $r=2$, then $\bar{K}\equiv0$ and, as a consequence, $\phi$ is constant. Therefore, $k_g\equiv0$, which is not possible since $\p\Omega$ is compact. So $r=1$ and the result follows.

\medskip

\textbf{Case II: $\Omega$ is compact unstable.} In this case, the result trivially holds for $r=1,2,3$. We only need to prove that $\Ind(\Omega)\ge\frac{r-2}{2}$ for $r\ge4$. Suppose, for contradiction, that $2n<r-2$, where $n=\Ind(\Omega)\ge1$.

Let $\lambda_1\le\cdots\le\lambda_n<0$ be the negative eigenvalues of $\Delta$ with corresponding linearly independent eigenfunctions $\phi_1,\ldots,\phi_n\in C^\infty(\Omega)$. We may assume that $\phi_1,\ldots,\phi_n$ are $L^2$-orthonormal.

Now define
\begin{align*}
\Phi:L^2\H_T^1(\Omega)\to\R^{2n},\quad\Phi(\xi)=\bigg(\int_\Omega\phi_i\langle\xi,E_j\rangle\bigg)_{\substack{i=1,\ldots,n\\ j=1,2}},
\end{align*}
where $\{E_1,E_2\}$ is the standard basis of $\R^2$. Since $\dim L^2\H_T^1(\Omega)=r-1>2n+1$, it follows from the rank-nullity theorem that $\dim\Ker\Phi\ge2$.

Given $\xi\in\Ker\Phi$, since $\langle\xi,E_j\rangle\in V^\perp$, we have
\begin{align*}
0\le\lambda_{n+1}\int_\Omega\langle\xi,E_j\rangle^2\le\int_\Omega|\nabla\langle\xi,E_j\rangle|^2+\int_{\p\Omega}k_g\langle\xi,E_j\rangle^2,\quad j=1,2.
\end{align*}
Summing over $j$, we obtain
\begin{align*}
0\le\lambda_{n+1}\int_\Omega|\xi|^2\le\int_\Omega|\nabla\xi|^2+\int_{\p\Omega}k_g|\xi|^2.
\end{align*}

On the other hand, it is not difficult to see that
\begin{align*}
\frac{1}{2}\Delta|\xi|^2=|\nabla\xi|^2+\langle\Delta\xi,\xi\rangle=|\nabla\xi|^2.
\end{align*}
Therefore,
\begin{align*}
0\le\lambda_{n+1}\int_\Omega|\xi|^2\le\int_{\p\Omega}(\langle\nabla_\nu\xi,\xi\rangle+k_g|\xi|^2).
\end{align*}

Now observe that, since $\langle\nabla_{\p_x}\xi,\p_y\rangle=\langle\nabla_{\p_y}\xi,\p_x\rangle$, then
\begin{align*}
\langle\nabla_\nu\xi,\xi\rangle=\langle\nabla_\xi\xi,\nu\rangle=-\langle\xi,\nabla_\xi\nu\rangle=-k_g|\xi|^2.
\end{align*}
Hence,
\begin{align*}
0\le\lambda_{n+1}\int_\Omega|\xi|^2\le\int_{\p\Omega}(\langle\nabla_\nu\xi,\xi\rangle+k_g|\xi|^2)=0.
\end{align*}

This forces each $\langle\xi,E_j\rangle$ to be an eigenfunction of $\Delta$, with Robin boundary condition $B\langle\xi,E_j\rangle=0$, associated with the eigenvalue $\lambda_{n+1}=0$. In particular, $\nabla_\nu\xi+k_g\xi=0$.

Denote $\xi=(a,b)$ and define $f:\Omega\subset\C\to\C$ by $f=a-ib$. It is not difficult to see that $\xi$ being harmonic is equivalent to $f$ being holomorphic.

Now fix a connected component $\Gamma$ of $\p\Omega$. Let $z(s)$ be an orientation-preserving arc-length parametrization of $\Gamma$ and denote $\tau(s)=z'(s)=e^{i\theta(s)}$, where $\theta(s)$ is the angle from $E_1$ to $z'(s)$. Straightforward computations give
\begin{align*}
f\tau=\langle\xi,\tau\rangle+i\langle\xi,\nu\rangle=\langle\xi,\tau\rangle.
\end{align*}

\textbf{Claim 1: $f\tau$ is constant on $\Gamma$.} Observe first that $\nabla_\nu\xi+k_g\xi=0$ is equivalent to $\p_\nu f+k_gf=0$. Therefore,
\begin{align*}
\frac{d}{ds}(fe^{i\theta})=(\p_\tau f+i\theta'f)e^{i\theta}=(\p_\nu f+k_gf)ie^{i\theta}=0,
\end{align*}
where we used that $\tau=i\nu$ (in particular, $\p_\tau f=i\p_\nu f$) and $\theta'=k_g$.

Thus, by Claim 1, we may define a linear map $\Psi:\Ker\Phi\to\R$ by letting $\Psi(\xi)$ be the constant value of $\langle\xi,\tau\rangle$ along $\Gamma$. Since a harmonic vector field (or a holomorphic function) is trivial when vanishing along a boundary component (see e.g.\,\cite[Theorem~3.4.4]{Schwarz}), we have $\Ker(\Psi)=\{0\}$, which is a contradiction because $\dim\Ker{\Phi}\ge2$, as we wanted to prove.

\medskip

\textbf{Case III: $\Omega$ is noncompact stable.} We claim that this case is empty. In fact, assume that $\Omega$ is stable. It can be proved that, in this case, there exists a positive harmonic function $\phi:\Omega\to\R$ satisfying the Robin boundary condition $B\phi=0$ (see \cite[Proposition~4.1]{HongSaturnino}; see also \cite[Proposition~2.1]{MazetMendes}). 

The metric $\bar{g}=\phi^2\delta$ has Gaussian curvature $\bar{K}\ge0$ and geodesic curvature $\bar{k}_g=0$:
\begin{align*}
\bar{K}\phi^2=-\Delta\log\phi=\frac{|\nabla\phi|^2}{\phi^2},\quad \bar{k}_g\phi=k_g+\p_\nu\log\phi=0.
\end{align*}
Also, the same arguments as in the proof of Theorem~1.4 of \cite{HongSaturnino} give that $\bar{g}$ is complete. Therefore, the Cohn-Vossen theorem for surfaces with boundary (see e.g.\,\cite[Theorem~2.2.1]{ShiohamaShioyaTanaka}) implies that 
\begin{align*}
0\le\int_\Omega\bar{K}+\int_{\p\Omega}\bar{k}_g\le2\pi\chi(\Omega)=2\pi(1-r).
\end{align*} 
Therefore, $r=1$ and $\bar{K}\equiv0$. This forces $\phi$ to be constant and, as a consequence, $k_g\equiv0$, which is not possible.

\medskip

\textbf{Case IV: $\Omega$ is noncompact unstable.} The result trivially holds in this case for $r=1,2,3$. As in the proof of Case II, assume for contradiction that $2n<r-2$.

Similarly to the compact case, there exists a linear subspace $V\subset L^2(\Omega)$ of $n$ dimensions, whose $L^2$-orthonormal basis consists of smooth eigenfunctions $\phi_1,\ldots,\phi_n$ of $\Delta$ associated with negative eigenvalues $\lambda_1\le\cdots\le\lambda_n<0$:
\begin{align*}
\begin{cases}
\Delta\phi_i+\lambda_i\phi_i=0&\mbox{in}\quad\Omega,\\
\p_\nu\phi_i+k_g\phi_i=0&\mbox{on}\quad\p\Omega.
\end{cases}
\end{align*}
Furthermore, 
\begin{align}\label{eq.aux.3}
Q(u,u)=\int_\Omega|\nabla u|^2+\int_{\p\Omega}k_gu^2\ge0
\end{align}
for every $u\in C_c^\infty(\Omega)\cap V^\perp$ (see \cite[Proposition~4.9]{HongSaturnino}).

As before, define the linear map
\begin{align*}
\Phi:L^2\H_T^1(\Omega)\to\R^{2n},\quad\Phi(\xi)=\bigg(\int_\Omega\phi_i\langle\xi,E_j\rangle\bigg)_{\substack{i=1,\ldots,n\\ j=1,2}},
\end{align*}
and observe that, by the rank-nullity theorem, we have $\dim\Ker(\Phi)\ge2$.

Now let $u_R:\R^2\to\R$, $R>0$, be a family of smooth functions satisfying the following conditions:
\begin{itemize}
\item $0\le u_R\le1$;
\item $u_R\equiv1$ on $D_R=\{x\in\R^2;|x|\le R\}$;
\item $u_R\equiv0$ on $\R^2\setminus D_{2R}$;
\item $|\nabla u_R|\le\frac{C}{R}$ for some constant $C>0$ independent of $R$;
\item $|\Delta u_R|\le\frac{C}{R^2}$.
\end{itemize}

For the existence of $u_R$, let $\psi:[0,+\infty)\to\R$ be any smooth function such that $0\le\psi\le 1$, $\psi\equiv1$ on $[0,1]$, $\psi\equiv0$ on $(2,+\infty)$, and define $u_R(x)=\psi(|x|/R)$. Since $\psi'$ and $\psi''$ have compact support, straightforward computations show that $u_R$ has the desired properties.

Let $X=(x_1,x_2):\Omega\to\R^2$ be a smooth vector field with compact support. Taking $R>0$ sufficiently large, we may assume that $\operatorname{supp}X\subset D_R$. Abusing notation, we write $X\in V^\perp$ to say that $x_1,x_2\in V^\perp$. With this notation in mind, $\Ker(\Phi)=L^2\H_T^1(\Omega)\cap V^\perp$. 

Fix $\xi\in\Ker(\Phi)$ and define $X_t=u_R(\xi+tX+\sum_{k=1}^n\phi_k\vec v_k)$, where $\vec v_1,\ldots,\vec v_n\in\R^2$ is the solution of 
\begin{align*}
\sum_{k=1}^n\bigg(\underbrace{\int_\Omega u_R\phi_i\phi_k}_{a_{ik}(R)}\bigg)v_{kj}=-\underbrace{\int_\Omega u_R\phi_i\langle\xi,E_j\rangle}_{b_{ij}(R)},\quad i=1,\ldots,n,
\end{align*}
with $\vec v_k=(v_{k1},v_{k2})$, which exists because, by the dominated convergence theorem, $a_{ik}(R)\to\delta_{ik}$ as $R\to+\infty$. Furthermore, $\vec v_k=\vec v_k(R)$ converges to the null vector as $R\to+\infty$, since $b_{ij}(R)\to0$. In other words, $u_R(\xi+\sum_{k=1}^n\phi_k\vec v_k)\in V^\perp$, which implies that $X_t\in V^\perp$ for every $t\in\R$, since $u_RX=X\in V^\perp$. 

Now, using the coordinates of $X_t$ as test-functions in \eqref{eq.aux.3}, we have:
\begin{align*}
0\le Q(X_t,X_t)&=Q(u_R\xi,u_R\xi)+2tQ(\xi,X)+2\sum_{k=1}^nQ(u_R\xi,u_R\phi_k\vec v_k)\\
&\qquad+t^2Q(X,X)+\sum_{k,\ell}^nQ(u_R\phi_k\vec v_k,u_R\phi_\ell\vec v_\ell),
\end{align*}
where we used that $\operatorname{supp}X\subset D_R$ and $Q(X,\phi_k\vec v_k)=0$, since $X\in V^\perp$. Here, for smooth vector fields $Y=(y_1,y_2)$ and $Z=(z_1,z_2)$ on $\Omega$, we define
\begin{align*}
Q(Y,Z)=\sum_{j=1}^2 Q(y_j,z_j),
\end{align*}
provided that at least one of them has compact support.

Thus, because the last inequality holds for every $t\in\R$, we deduce that
\begin{align*}
Q(\xi,X)^2\le Q(X,X)\Big(\underbrace{Q(u_R\xi,u_R\xi)}_{\rm (I)}&+2\sum_{k=1}^n\underbrace{Q(u_R\xi,u_R\phi_k\vec v_k)}_{\rm (II)}\\
&\qquad+\sum_{k,\ell=1}^n\underbrace{Q(u_R\phi_k\vec v_k,u_R\phi_\ell\vec v_\ell)}_{\rm (III)}\Big).
\end{align*}

So direct computations as in \cite{CavalcanteMendesdosSantos} or \cite{ChodoshMaximo-I} give that
\begin{align*}
{\rm (I)}&=\int_\Omega|\xi|^2|\nabla u_R|^2,\\
{\rm (II)}&=\lambda_k\int_\Omega u_R^2\langle\xi,\vec v_k\rangle\phi_k-\int_\Omega(u_R\Delta u_R)\langle\xi,\vec v_k\rangle\phi_k-2\int_\Omega u_R\langle\xi,\vec v_k\rangle\langle\nabla u_R,\nabla\phi_k\rangle,\\
{\rm (III)}&=\frac{\lambda_k+\lambda_\ell}{2}\langle\vec v_k,\vec v_\ell\rangle\int_\Omega u_R^2\phi_k\phi_\ell+\langle\vec v_k,\vec v_\ell\rangle\int_\Omega|\nabla u_R|^2\phi_k\phi_\ell.
\end{align*} 
Then, applying the dominated convergence theorem to each of the above integrals together with the fact that $\vec v_k(R)\to0$ as $R\to+\infty$, and using the properties of $u_R$, we conclude that (I), (II) and (III) go to zero as $R\to+\infty$. Remember that $|\xi|,\phi_k\in L^2(\Omega)$ and, since $\phi_k$ is an eigenfunction associated with $\lambda_k$,
\begin{align*}
\int_\Omega|\nabla\phi_k|^2=\lambda_k\int_\Omega\phi_k^2-\int_{\p\Omega}k_g\phi_k^2<+\infty.
\end{align*}

We have shown that
\begin{align}\label{eq.aux.1}
\int_{\p\Omega}\langle X,B\xi\rangle=-\int_\Omega\langle X,\Delta\xi\rangle+\int_{\p\Omega}\langle X,B\xi\rangle=Q(X,\xi)=0
\end{align}
for every smooth vector field $X\in V^\perp$ with compact support, where we denote $B\xi=\nabla_\nu\xi+k_g\xi$. Taking $X=u_RY$ and observing that $u_R\equiv1$ on $\p\Omega$ for $R>0$ sufficiently large, we see that \eqref{eq.aux.1} holds for every smooth vector field $Y\in V^\perp$.

On the other hand, observe that for each $\vec a\in\R^2$,
\begin{align*}
\int_{\p\Omega}\langle\phi_k\vec a,B\xi\rangle&=\int_{\p\Omega}\langle u_R\phi_k\vec a,B\xi\rangle\\
&=-\int_\Omega\langle\xi,\Delta(u_R\phi_k\vec a)\rangle+\int_{\p\Omega}\langle\xi,B(u_R\phi_k\vec a)\rangle\\
&=\lambda_k\int_\Omega u_R\langle\xi,\vec a\rangle\phi_k-\int_\Omega (\Delta u_R)\langle\xi,\vec a\rangle\phi_k-2\int_\Omega\langle\xi,\vec a\rangle\langle\nabla u_R,\nabla\phi_k\rangle,
\end{align*}
where we used that $\Delta\xi=0$, $u_R\equiv1$ on $\p\Omega$, $\Delta\phi_k=-\lambda_k\phi_k$ in $\Omega$, and $B(\phi_k)=0$ along $\p\Omega$. Therefore, since $u_R\to1$, $\Delta u_R\to0$, and $|\nabla u_R|\to0$ as $R\to+\infty$, and since $\langle\xi,\vec a\rangle\in V^\perp$, it follows from the dominated convergence theorem that the right-hand side converges to zero. Hence,
\begin{align}\label{eq.aux.2}
\int_{\p\Omega}\langle\phi_k\vec a,B\xi\rangle=0,\quad\forall\vec a\in\R^2,\quad \forall k=1,\ldots,n.
\end{align}

Identity \eqref{eq.aux.1} applied to $Y$, together with \eqref{eq.aux.2}, implies that $B\xi=\nabla_\nu\xi+k_g\xi=0$ on $\p\Omega$. Therefore, arguing exactly as at the end of the proof of Case II, we obtain a contradiction. This finishes the proof of Theorem~\ref{thm.main}.

\section{Examples}

\subsection{The unit disk.}

We want to calculate the first eigenvalue of the unit disk
\begin{align*}
\D=\{z\in\R^2;|z|\le 1\}.
\end{align*}

Making the separation of variables $\phi=R(r)e^{im\theta}$, $m\in\Z$, the problem \eqref{problema} is equivalent to
\begin{align*}
\begin{cases}
r^2R''+rR'+(\lambda r^2-m^2)R=0&\mbox{in}\quad[0,1],\\
R'(1)+R(1)=0.
\end{cases}
\end{align*}

For $\lambda>0$, writing $\lambda=\mu^2$, the regular solutions are given by the Bessel functions of the first kind,
\begin{align*}
R(r)=J_m(\mu r).
\end{align*}
The boundary condition implies
\begin{align}\label{eq.aux.5}
\mu J_m'(\mu)+J_m(\mu)=0.
\end{align}

Therefore, the eigenvalues of $\Delta$ with Robin boundary condition $B\phi=0$ on the disk are given by
\begin{align*}
\lambda_{m,k}=\mu_{m,k}^2,
\end{align*}
where $\mu_{m,k}$ denotes the $k$-th positive root of \eqref{eq.aux.5}. For $m\ge 1$, the functions
\begin{align*}
J_m(\mu_{m,k}r)\cos(m\theta)
\quad\text{and}\quad
J_m(\mu_{m,k}r)\sin(m\theta)
\end{align*}
are linearly independent eigenfunctions associated with $\lambda_{m,k}$, so these eigenvalues have multiplicity at least two. Hence the first eigenvalue must correspond to the radial mode $m=0$. Numerical computations give
\begin{align*}
\mu_{0,1}\approx1.2558,
\end{align*}
and therefore
\begin{align*}
\lambda_1=\lambda_{0,1}=\mu_{0,1}^2\approx1.5770.
\end{align*}

\subsection{The annulus.}

Consider the annulus
\begin{align*}
\mathbb{A}=\{z\in\R^2;1\le |z|\le 2\}.
\end{align*}

Using the Fourier decomposition
\begin{align*}
u(r,\theta)=\sum_{m\in\Z}R_m(r)e^{im\theta},
\end{align*}
the quadratic form
\begin{align*}
Q(u)=\int_\A|\nabla u|^2+\int_{\p\A}k_g u^2=\int_0^{2\pi}\int_1^2|\nabla u|^2rdrd\theta+\frac12\int_{S^1(2)}u^2ds-\int_{S^1(1)}u^2ds
\end{align*}
decomposes as
\begin{align*}
Q(u)=2\pi\sum_{m\in\Z}Q_m(R_m),
\end{align*}
where
\begin{align*}
Q_m(R)=\int_1^2\bigg(r(R')^2+\frac{m^2}{r}R^2\bigg)dr+R(2)^2-R(1)^2.
\end{align*}

Observe that $Q_{-m}=Q_m$ and
\begin{align*}
Q_m(R)-Q_n(R)=(m^2-n^2)\int_1^2\frac{R^2}{r}dr>0
\end{align*}
for every $R\not\equiv0$ whenever $m>n\ge0$. Therefore, it suffices to analyze the modes $m=0$ and $m=1$.

We first show that $Q_1$ is nonnegative. Indeed,
\begin{align*}
r\bigg(R'+\frac{R}{r}\bigg)^2=r(R')^2+\frac{R^2}{r}+(R^2)'.
\end{align*}
So, integrating over $[1,2]$, we obtain
\begin{align*}
Q_1(R)=\int_1^2r\bigg(R'+\frac{R}{r}\bigg)^2dr\ge0.
\end{align*}

Moreover, $Q_1(R)=0$ if and only if $R(r)=A/r$, $A\in\R$. This gives rise to the two functions
\begin{align*}
\frac{\cos\theta}{r}=\frac{x}{x^2+y^2},\quad\frac{\sin\theta}{r}=\frac{y}{x^2+y^2},
\end{align*}
which belong to $\Ker(Q)$.

We now analyze the quadratic form
\begin{align*}
Q_0(R)=\int_1^2 r(R')^2dr+R(2)^2-R(1)^2.
\end{align*}

Define $v(r)=R(r)-R(1)$, so that $R(r)=R(1)+v(r)$ and $v(1)=0$. This yields the decomposition
\begin{align*}
C^\infty([1,2])=\operatorname{span}\{1\}\oplus H,\quad H=\{v\in C^\infty([1,2]);v(1)=0\}.
\end{align*}

Now observe that
\begin{align*}
Q_0(v)=\int_1^2 r(v')^2dr+v(2)^2>0,\quad\forall v\in H\setminus\{0\}.
\end{align*}

Since $H$ has codimension one and $Q_0$ is positive on $H$, any subspace on which $Q_0$ is negative definite must intersect $H$ trivially. Therefore, $\Ind(Q_0)\le1$. On the other hand, taking $R(r)=3-r$, a direct computation gives $Q_0(R)=-\frac32<0$. Hence $\Ind(Q_0)=1$.

In summary, the calculations above show that $\lambda_1<0$, $\lambda_2=\lambda_3=0$, and $\lambda_k>0$ for $k\ge4$. Moreover, the eigenspace corresponding to the zero eigenvalue is precisely
\begin{align*}
\operatorname{span}\bigg\{\frac{x}{x^2+y^2},\frac{y}{x^2+y^2}\bigg\}.
\end{align*}

\subsection{The exterior of the disk} In this final example, we consider the planar domain outside the open unit disk
\begin{align*}
\E=\{z\in\R^2:|z|\ge1\}.
\end{align*}
We compute the first eigenvalue and the index of $\E$.

As in the case of the disk, after separation of variables, problem \eqref{problema} reduces to
\begin{align}\label{eq.plane}
\begin{cases}
r^2R''+rR'+(\lambda r^2-m^2)R=0&\mbox{in}\quad[1,+\infty),\\
R'(1)+R(1)=0.
\end{cases}
\end{align}

We are interested in negative eigenvalues $\lambda=-\alpha^2<0$, where $\alpha>0$. In this case, the differential equation in \eqref{eq.plane} can be written as
\begin{align*}
r^2R''+rR'-(\alpha^2r^2+m^2)R=0\qmbox{in}[1,+\infty),
\end{align*}
whose general solution is given in terms of the modified Bessel functions of the first and second kind:
\begin{align*}
R(r)=AI_m(\alpha r)+BK_m(\alpha r),
\end{align*}
where $I_m$ has exponential growth and $K_m$ has exponential decay as $r\to+\infty$. Since we are looking for eigenfunctions in $L^2(\E)$, we must have $A=0$.

The boundary condition in \eqref{eq.plane} becomes
\begin{align*}
\alpha K_m'(\alpha)+K_m(\alpha)=0.
\end{align*}
Using the identity
\begin{align*}
K_m'(\alpha)=-K_{m-1}(\alpha)-\frac{m}{\alpha}K_m(\alpha),
\end{align*}
we obtain
\begin{align*}
\alpha K_m'(\alpha)+K_m(\alpha)=-\alpha K_{m-1}(\alpha)-(m-1)K_m(\alpha).
\end{align*}
Since $K_m(\alpha)>0$, it follows that
\begin{align*}
\alpha K_m'(\alpha)+K_m(\alpha)<0
\end{align*}
for every $m\ge1$ and $\alpha>0$. Because $K_{-m}=K_m$, this shows that there are no negative eigenvalues corresponding to the modes $m\neq0$.

For $m=0$, we seek positive roots of
\begin{align*}
\alpha K_0'(\alpha)+K_0(\alpha)=0.
\end{align*}
Numerical computation shows that this equation has a unique positive root,
\begin{align*}
\alpha_{0,1}\approx0.5950.
\end{align*}
Therefore, $\Ind(\E)=1$ and
\begin{align*}
\lambda_1=-\alpha_{0,1}^2\approx-0.3540.
\end{align*}

A similar approach to that of the annulus also shows that $\Ind(\E)=1$.

\begin{remark}
The planar domains $\D$, $\A$, and $\E$ can be realized as free boundary\linebreak
minimal surfaces in suitable domains of $\R^3$ whose boundaries are formed by catenoids. More precisely, the corresponding ambient spaces are the interior of a catenoid, the region bounded by two concentric catenoids, and the exterior of a catenoid, respectively.
\end{remark}

\begin{remark}
It would be interesting to determine whether equality can occur in the index estimate of Theorem~\ref{thm.main}. Although the examples above are useful for illustrating the computation of the index and the first eigenvalue in concrete cases, they do not actually help in addressing this question.
\end{remark}

\begin{acknowledgments}
This study was financed in part by the Conselho Nacional de Desenvolvimento Científico e Tecnológico (CNPq, Grants No.~309867/2023-1 and 445723/2025-4).
\end{acknowledgments}

\subsection*{Data availability}
No data were generated or analyzed in this study.

\bibliographystyle{amsplain}
\bibliography{bibliography.bib}

\end{document}